\documentclass[11pt]{amsart}
\usepackage{texmac}

\bbsymbols{}{A,P,G}
\bbsymbols{b}{H}
\calsymbols{}{E,L,T,B}
\calsymbols{c}{C,A,P,O,K,G,F}
\def\H{H_{24}}
\begin{document}

\vskip-0mm\llap{.\hskip10cm}

\title{Surjectivity of mod $2^n$ representations of elliptic curves}
\author{Tim and Vladimir Dokchitser}
\address{Robinson College, Cambridge CB3 9AN, United Kingdom}
\email{t.dokchitser@dpmms.cam.ac.uk}
\address{Emmanuel College, Cambridge CB2 3AP, United Kingdom}
\email{v.dokchitser@dpmms.cam.ac.uk}
\date{April 25, 2011}

\maketitle

Associated to an elliptic curve $E/\Q$ and a prime $l$ are representations
\beq
  &\qquad&\bar\rho_{l^n}: &\Gal(\bar\Q/\Q)\lar \Aut E[l^n] &\iso \GL_2(\Z/l^n\Z)\\[5pt]
\noalign{\hskip-1em and their inverse limit, the $l$-adic representation}
  &\vphantom{\int^{X^{A^A}}}\qquad&\rho_l: &\Gal(\bar\Q/\Q)\lar\Aut T_l E &\iso \GL_2(\Z_l).
\eeq
As explained by Serre (\cite{SerA} IV, 3.4), for $l\ge 5$
the group $\SL_2(\Z_l)$ has no proper closed subgroups that surject
onto $\SL_2(\F_l)$, so
\beq
  &\qquad&
  \bar\rho_l\text{ surjective} & \implies & \rho_l\text{ surjective}
  &&\text{for }l\ge 5.\\[5pt]
\noalign{\hskip-1em The condition $l\ge 5$ is necessary, and}
  \vphantom{\int^{X^{A^A}}}
  &\qquad&
  \bar\rho_{l^n}\text{ surjective} & \>\>\not\!\!\!\implies & \bar\rho_{l^{n+1}}\text{ surjective}
  &&\text{for }l^n=2,3,4.
\eeq
Elliptic curves with surjective mod 3 but not mod 9 representation have been
classified by Elkies \cite{Elk}, and the purpose of this note is to do this
in the
`2 not 4' and `4 not 8' cases as well:

\begin{theorem*}
Let $E: y^2=x^3+ax+b$ be an elliptic curve over $\Q$ with discriminant
$\Delta=-16(4a^3+27b^2)$ and $j$-invariant $j=-1728(4a)^3/\Delta$.~Then
\begin{enumerate}
\item
$\bar\rho_2$ is surjective $\iff$ $x^3+ax+b$ is irreducible
and $\Delta\notin\Q^{\times 2}$.
\item
$\bar\rho_4$ is surjective $\iff$ $\bar\rho_2$ is surjective,
$\Delta\notin-1\cdot\Q^{\times 2}$ and
$
  j \ne -4t^3(t+8)
$
for any $t\in\Q$.
\item
$\bar\rho_8$ is surjective $\iff$ $\bar\rho_4$ is surjective and
$\Delta\notin\pm 2\cdot\Q^{\times 2}$.
\end{enumerate}
\end{theorem*}

\begin{proof}
The $x$-coordinates of the three non-trivial
2-torsion points are the roots of $x^3+ax+b$ and their $y$-coordinates are 0.
So $\bar\rho_2$
surjects onto $\GL_2(\F_2)\iso S_3$ if and only if this cubic is
irreducible and its discriminant $\Delta/16$ is not a square.
Note that this proves (1) and that $\Q(E[2])$ contains $\Q(\sqrt\Delta)$.

Now recall that by the properties of the Weil pairing,
$\Q(E[n])\supset \Q(\zeta_n)$ and the corresponding map
$\Gal(\Q(E[n])/\Q)\surjects(\Z/n\Z)^\times$ is simply the determinant.
In particular,
$$
  \Q(E[4])\supset\Q(\sqrt\Delta,\sqrt{-1})
    \quad\text{and}\quad
  \Q(E[8])\supset\Q(\sqrt\Delta,\sqrt{-1},\sqrt{2}).
$$
Incidentally, as there are elliptic curves whose 2-torsion defines an
$S_3$-exten\-sion of $\Q$ which is disjoint from $\Q(\zeta_8)$, e.g.
$y^2=x^3-2$, this shows that the canonical maps $(\hskip-0.1em\bmod\>2,\det)$
$$
  \GL_2(\Z/4\Z)\lar S_3\times(\Z/4\Z)^\times
    \quad\text{and}\quad
  \GL_2(\Z/8\Z)\lar S_3\times(\Z/8\Z)^\times
$$
are surjective (being already surjective on the subgroups
$\Im\bar\rho_4$ and $\Im\bar\rho_8$ for this elliptic curve).
Hence, if an elliptic curve $E/\Q$ has surjective $\bar\rho_4$,
then $\Q(E[2],\zeta_4)$ has degree $12$, so $\Q(\sqrt\Delta,\sqrt{-1})$
has degree 4. Similarly, if $\bar\rho_8$ is surjective, then
$\Q(\sqrt\Delta,\sqrt{-1},\sqrt{2})$ has degree 8.

\smallskip

(3)
If $\bar\rho_8$ is surjective, then so is $\bar\rho_4$, and as
$\Q(\sqrt\Delta)\ne \Q(\sqrt{\pm2})$, it follows that
$\Delta\notin\pm 2\cdot\Q^{\times 2}$.
Conversely, if $\bar\rho_4$ is surjective and $\Delta\notin\pm 2\cdot\Q^{\times2}$,
then $\Q(\sqrt\Delta,\sqrt{-1},\sqrt{2})$ is a
$C_2\times C_2\times C_2$-extension of $\Q$.
So $\Im\bar\rho_8$ surjects onto $\GL_2(\Z/4\Z)$ and onto $(\Z/8\Z)^\times$,
and possesses a $C_2\times C_2\times C_2$-quotient.
A computation shows that the
only such subgroup of $\GL_2(\Z/8\Z)$ is the full group itself.


\smallskip

(2) The argument is the same as for (3), except that in this case
$\GL_2(\Z/4\Z)$ does have a (unique up to conjugacy) proper subgroup
which surjects onto $\GL_2(\Z/2\Z)$ and onto $(\Z/4\Z)^\times$,
and has a $C_2\times C_2$-quotient.
This group has index 4, and is conjugate to
$\H=\langle\smallmatrix 0130, \smallmatrix 0111\rangle \iso C_3\rtimes D_8$.
The following lemma completes the proof.
\end{proof}

\begin{lemma*}
Let $E/\Q$ be the elliptic curve $y^2=x^3+ax+b$ with $b\ne 0$. The following conditions are equivalent:
\begin{enumerate}
\item $\Gal(\Q(E[4])/\Q)$ is conjugate to a subgroup of $\H$.
\item The polynomial
$$
  \qquad f(x)=x^4-4ax^3+6a^2x^2 + 4(7a^3+54b^2)x + (17a^4+108ab^2)
$$
has a rational root.
\item $j(E)=-4t^3(t+8)$ for some $t\in\Q$.
\end{enumerate}
\end{lemma*}

\begin{proof}
(1)$\iff$(2). Regard $a$ and $b$ as indeterminants and
$E$ as a curve over $K=\Q(a,b)$.
Consider the 4-torsion polynomial (cf. $\psi_4$ in \cite{Sil1} Exc. III.3.7)
$$
  \psi(x) = x^6 + 5 ax^4 + 20bx^3 - 5a^2x^2 - 4abx - (a^3+8b^2).
$$
Its roots $x_1,...,x_6\in\bar K$ are the $x$-coordinates of the
primitive 4-torsion points of~$E$.
Pick a basis $P=(x_1,y_1), Q=(x_2,y_2)$ for the 4-torsion, and
for each of the four left cosets $C$ of $\H$ in $\GL_2(\Z/4\Z)$
define
$$
   \theta_C = \sum_{g\in C} x(g P)x(g Q),
$$
where $x(\cdot)$ is the $x$-coordinate.
The $\theta_C$ are distinct, which can be checked by specialising
e.g. $a\mapsto 0$, $b \mapsto 1$, where the numbers are
$0, -6, 3\pm 3\sqrt{-3}$.

Now for $h\in\Gal(K(E[4])/K)\subset \GL_2(\Z/4\Z)$ we have
$h(\theta_C)=\theta_{hC}$, and hence the polynomial
$$
  \tilde f(x)=\prod_C (x-\theta_C)=\smash{\sum_{i=0}^4} \lambda_i x^i
$$
has coefficients in $K$.
The roots $x_i\in\bar K$ of $\psi$ are integral over $R=\Z[a,b]$,
and therefore so are the $\theta_C$.
Because $R$ is integrally closed in its
field of fractions $K$ and
the $\lambda_i$ are both in $K$ and are integral over~$R$,
they lie in~$R$.
Moreover, if we rescale $a\mapsto s^2a, b\mapsto s^3b$, then the roots
of $\psi(x)$ change to $x_i\mapsto sx_i$, and therefore
$\theta_C\mapsto s^2\theta_C$.
So $\lambda_i$ must have weight $2i$, in the sense that
$$
  \lambda_1\in \Z a, \quad \lambda_2\in \Z a^2, \quad \lambda_3\in \Z a^3+\Z b^2 \quad
  \text{and}\quad
  \lambda_4\in \Z a^4+\Z ab^2.
$$
Now we can compute $\tilde f$ numerically for a few specialisations $a,b\in\Z$ (using
complex uniformisation of torsion points and rounding the coefficients) to deduce the
exact formulae for the $\lambda_i$, and we find that $\tilde f(x)=f(4x)$.

Let us return to $a,b\in\Q$, constructing $\psi, \theta_C$ and $\tilde f$
in the same way. Note that
the discriminant of $f$ is
$3^6 b^2 \Delta_E^3$,
so it has no repeated roots and the $\theta_C$ are distinct.
Now $h: \theta_C\mapsto\theta_{hC}$ defines a
transitive action
of $\GL_2(\Z/4\Z)$ on these, and for $h\in\Gal(\Q(E[4])/\Q)$
it coincides with the Galois action $\theta_C\mapsto h(\theta_C)$.
Because all four $\theta_C$ are distinct,
the stabiliser of $\theta_{\H}$ is precisely $\H$, and
the stabilisers of the others are the conjugates of $\H$.
So $\Gal(\Q(E[4])/\Q)$ is conjugate to a subgroup of $\H$ if and only if
one of the $\theta_C$ is rational, equivalently if $f$ has a rational root.

(2)$\iff$(3).
First note that if $a=0$, equivalently $j=0$, then both conditions are satisfied
($f(0)=0$ and $j=0$). Suppose $a\ne 0$ and that
$f(x)$ has a rational root $r$. Then $u=r/a$ satisfies
$$
  u^4-4u^3+6u^2 + 4(7+54\frac{b^2}{a^3})u + (17+108\frac{b^2}{a^3}) = 0.
$$
Rewriting $\frac{b^2}{a^3}$ in terms of $j=\frac{1728(4a)^3}{16(4a^3+27b^2)}$
(namely $\frac{b^2}{a^3}=-\frac{4(j-1728)}{27j}$),
we find that $j=-27648\frac{2u+1}{(u-1)^4}$. Replacing $u$ by $\frac{12}t+1$ we get $j=-4t^3(t+8)$, as claimed.
Reversing the argument gives the other implication as well.
\end{proof}

\begin{remark*}
The elliptic curve
$$
  y^2 = x^3 - \frac{3t^2 + 24t}{t^2 - 4t + 12}x - \frac{2t^2 + 28t + 96}{t^2 - 4t + 12}
  \eqno{(t\ne -8)}
$$
has $j$-invariant $-4t^3(t+8)$. So for every curve of this form
the polynomial $f(x)$ has a rational root and, conversely,
every elliptic curve over $\Q$ for which $f(x)$ has a rational root
is a twist of a curve in the family.
\end{remark*}

\begin{acknowledgements}
We would like to thank Matt Greenberg and Gagan Sekhon for asking the
questions that led to this work. The first author is supported by a
Royal Society University Research Fellowship.
\end{acknowledgements}

\end{document}